\documentclass[11pt]{article}
\usepackage{a4}
\usepackage{amssymb}
\usepackage{graphicx}

\def\PP{{\mathbb P}}
\def\RR{{\mathbb R}}
\def\CC{{\mathbb C}}

\def\ZZ{{\mathbb Z}}
\def\wt#1{\widetilde{#1}}
\def\la{\lambda}
\def\ep{\varepsilon}

\newtheorem{lemma}{Lemma}

\newtheorem{proposition}[lemma]{Proposition}

\newcommand{\cn}{\mathop{\rm :}}
\def\roep #1.{\medbreak\noindent{\em #1\/}.\enspace}
\def\endroep{\par\medbreak}
\def\qedsymbol{$\Box$}
\newcommand\qed{\relax\ifmmode\qedsymbol\else
  {\unskip\nobreak\hfil\penalty50\hskip1em\null\nobreak\hfil\qedsymbol
  \parfillskip=0pt\finalhyphendemerits=0\endgraf\endroep}\fi}

\long\def\comment#1\endcomment{}

\begin{document}

\title{\bf Sextic surfaces with ten triple points}
\author{Jan Stevens
       \thanks{Partially supported by the 
               Swedish Research Council (VR)}}
\date{}
\maketitle

\begin{abstract}
  All families of sextic surfaces with the maximal number of isolated
  triple points are found.
\end{abstract}

Surfaces in $\PP^3(\CC)$ with isolated ordinary triple points have
been studied  in  \cite{EPS}. 
The results are most complete for degree six. 
A sextic surface can have at most ten triple points, and such surfaces
exist. For up to
nine triple points \cite{EPS} contains a complete classification.
In  this note I achieve the same for ten triple points.

The study of sextics with nine triple points is easier, because
they do lie on a quadric $Q$, which is not the case for ten points. 
Given such a sextic with equation $F$ the 
general element of the pencil $\alpha F + \beta Q^3$ is again a sextic
with nine isolated triple points. It turns out that such a pencil
also contains reducible surfaces, which are much easier to construct.
The same argument shows that a sextic with ten triple points is a
degeneration of one with nine (simply choose a quadric
through nine of the ten points).

Therefore one can look for sextics with ten triple points in each of the
five families given in \cite{EPS}. It suffices to consider only those
two,which have a rather nice description. 
The  one-parameter family of examples  \cite{EPS} was found in
the first family by imposing extra symmetry. The surfaces in the other family 
have the simplest equations of all. Nevertheless I could not find
a single solution, because I was looking at the wrong place: as explained
below, I made an unwarranted general position assumption. Different
families of sextics are connected by Cremona transformations. By
transforming the known example I found the right
assumptions.  The equations for a tenth triple point in the family
become very simple, as I had hoped all the time.   
They describe a three-dimensional family. Knowing the dimension
then helped to find all solutions in the other family.

To solve the equations I use the
computer algebra system  {\tt Singular} \cite{Singular}.
The equations for the tenth point come from the ten second partial
derivatives of the defining function. 
The families with nine triple points depend on seven or eight
moduli. The unknown position of the tenth singular point adds three
more variables.  
One gets a very complicated system, which only can be
attacked by using the special structure of the equations.

The main result is that there are four different families of
sextics with ten triple points, each depending on three moduli.
They are distinguished by the number of $(-1)$-conics, which ranges
from two to five.


\section{Nine triple points}
The clue to the classification of  sextics with many triple points is
the study of exceptional curves of the first kind on the minimal
resolution. 
Let $X$ be a sextic with isolated triple points and $\wt X$ its
minimal resolution.
Whenever the canonical divisor $K_{\wt X}$ is effective,
any exceptional curve of the first kind $E$ is automatically
a  component, as $K_{\wt X}\cdot E = -1$.
Therefore $E$ comes from a rational curve on $X$
which is contained in the base locus of the system of quadrics 
through the triple points.
Assume that $X$ has nine  triple points $P_1,\ldots,P_9$. Let
$Q$ be the unique (irreducible) canonical quadric surface
and let $K=Q\cdot X$ be the adjoint curve.
The resolution $\wt X$ has exactly three disjoint $(-1)$-curves 
$C_1,C_2,C_3$ of
degrees $c_1,c_2,c_3\in\{2,4,5,6,7,8\}$ which are components of
$K$.
There are two 
possibilities: either $C_1+C_2+C_3=K$ or not.
In the first case $\wt X$ is a $K3$ surface blown up in
three points. By \cite{EPS}, Prop.~4.10,
there are up to permutation three choices for the degrees: 
$$
    (c_1,c_2,c_3)\in\{(2,2,8),(2,4,6),(4,4,4)\}.
$$
In the second case we end up with an effective canonical divisor
after blowing down $C_1$, $C_2$ and $C_3$. Now $\wt X$ is the blowup of
a minimal properly elliptic surface in three points and
by \cite{EPS}, Prop.~4.9, up  to permutation
$$
    (c_1,c_2,c_3)\in\{(2,2,2),(2,2,4)\}.
$$
In all cases
the curves $C_i$ of degree $c_i$ can be constructed as complete intersection
of $Q$ and a surface of degree $c_i/2$. In particular, if $c_i=2$ we have
five points on a conic in a plane. Such a conic will be called
a $(-1)$-conic. We call the triple $(c_1,c_2,c_3)$ the {\sl type} of
the surface.  

\begin{figure}
\center
\includegraphics{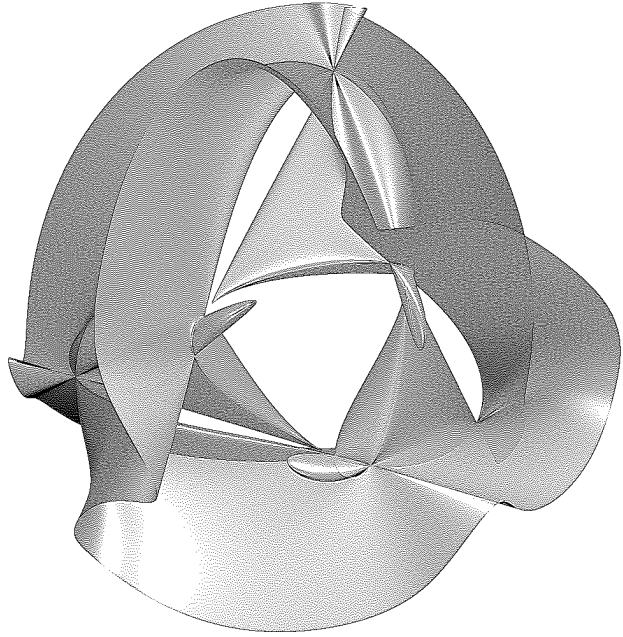}
\caption{A surface of type $(4,4,4)$ with nine triple points}
\label{surface}
\end{figure}

For every $(c_1,c_2,c_3)\in\{(2,2,8),(2,4,6),(4,4,4)\}$ there
exists a seven parameter family of sextic surfaces with
nine triple points and three $(-1)$-curves of degrees
$c_1$, $c_2$ and $c_3$ (\cite{EPS}, Thm.~4.13). 
Moreover $X$ occurs in a pencil of the form
$$
      \begin{array}{r@{\quad}l}
      \alpha\,  K_1K_2K_3+ \beta\, Q^3=0\;, &  
                              \mathrm{if}\ (c_1,c_2,c_3)=(4,4,4)\;,\\[1ex]
      \alpha\,  L_1K_2C_3+ \beta\, Q^3=0\;, & 
                              \mathrm{if}\ (c_1,c_2,c_3)=(2,4,6)\;,\\[1ex]
      \alpha\,  L_1L_2Q_3+\beta\,  Q^3=0\;, & 
                              \mathrm{if}\ (c_1,c_2,c_3)=(2,2,8)\;.
      \end{array}
$$
Here $Q$ is the unique canonical surface, 
$L_i$ stands for a linear form, $K_i$ for a singular quadric,
$C_3$ defines a four nodal cubic and $Q_3$  a
quartic surface with a triple point and six double points.
The multiplicities of the three surfaces in the nine singular points
are displayed in  Table \ref{table:multiK3}.
Note that we do not distinguish between a surface
and the form defining it, which we also call its equation.
\begin{table}
    \centering
    \begin{tabular}{|l|c||c|c|c|c|c|c|c|c|c|}\hline
       type & surface
       &$P_1$&$P_2$&$P_3$&$P_4$&$P_5$&$P_6$&$P_7$&$P_8$&$P_9$\\\hline\hline
        & $K_1$ & 0 & 2 & 1 & 1 & 1 & 1 & 1 & 1 & 1 \\\cline{2-11}
       {$(4,4,4)$}
        & $K_2$ & 1 & 0 & 2 & 1 & 1 & 1 & 1 & 1 & 1 \\\cline{2-11}
        & $K_3$ & 2 & 1 & 0 & 1 & 1 & 1 & 1 & 1 & 1 \\\hline\hline
        & $L_1$ & 1 & 0 & 0 & 0 & 0 & 1 & 1 & 1 & 1 \\\cline{2-11}
      {$(2,4,6)$} 
         & $K_2$ & 0 & 2 & 1 & 1 & 1 & 1 & 1 & 1 & 1 \\\cline{2-11}
        & $C_3$ & 2 & 1 & 2 & 2 & 2 & 1 & 1 & 1 & 1 \\\hline\hline
        & $L_1$ & 1 & 0 & 0 & 0 & 0 & 1 & 1 & 1 & 1 \\\cline{2-11}
        {$(2,2,8)$} 
         & $L_2$ & 0 & 1 & 1 & 1 & 0 & 0 & 0 & 1 & 1 \\\cline{2-11}
        & $Q_3$ & 2 & 2 & 2 & 2 & 3 & 2 & 2 & 1 & 1 \\\hline\hline
    \end{tabular}
    \caption{multiplicities at the singular points in the $K3$-case}
    \label{table:multiK3}
\end{table}
Figure \ref{surface} shows a surface of type $(4,4,4)$. The picture
was made with Stephan Endra\ss' program \texttt{surf} \cite{surf}.

\begin{figure}
\center
\includegraphics{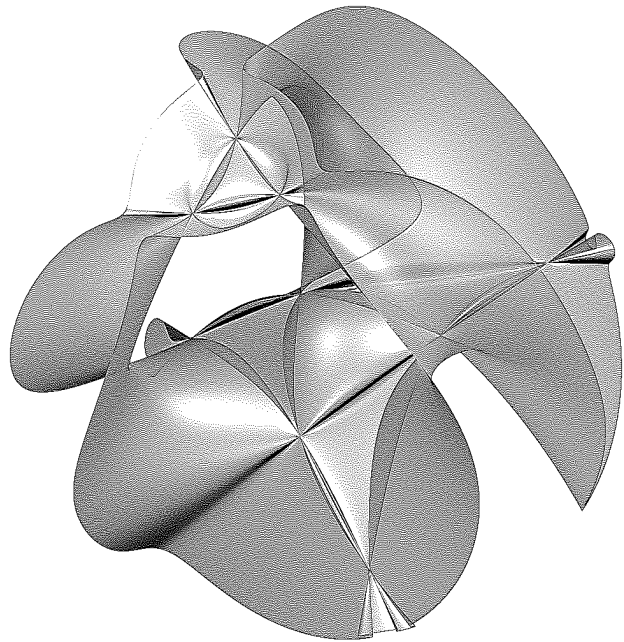}
\caption{A surface of type $(2,2,2)$ with nine triple points}
\label{surfacf}
\end{figure}

For every $(c_1,c_2,c_3)\in\{(2,2,2),(2,2,4)\}$ there
exists an eight parameter family of sextic surfaces with
nine triple points and three $(-1)$-curves of degrees
$c_1$, $c_2$ and $c_3$ (\cite{EPS}, Thm.~4.14). 
Moreover $X$ occurs in a web of the form
$$
  \begin{array}{r@{}c@{}c@{\quad}l}
    \alpha\, L_1L_2L_3 C  + {}&\beta\, L_1L_2L_3HQ& {}+ \gamma\, Q^3=0\;,
           & \mathrm{if}\  (c_1,c_2,c_3)=(2,2,2)\;,\\[1ex]
    \alpha \,L_1L_2K_3 Q'  + {}&\beta\, L_1L_2K_3Q & {}+\gamma\, Q^3=0\;,
          & \mathrm{if}\  (c_1,c_2,c_3)=(2,2,4)\;.
  \end{array}
$$
%
Again $L_i$ stands for a linear form. 
In the case $(2,2,2)$ the plane $H$ passes through
the three triple points not lying on  the double lines of $L_1L_2L_3$.
The  reducible cubic $HQ$ is an element of the pencil
of cubics through all points with double points
in $P_7$, $P_8$ and $P_9$, and $C$ is another such cubic. 
In the case $(2,2,4)$ the surface $K_3$ is a quadric cone and
$Q'$ is a smooth quadric not passing through $P_6$.
The multiplicities in the nine triple points of the surfaces giving
$(-1)$-curves  are displayed in Table~%
\ref{table:multelliptic}. A surface of type $(2,2,2)$
is shown in Figure \ref{surfacf}.
\begin{table}
    \centering
    \begin{tabular}{|l|c||c|c|c|c|c|c|c|c|c|}\hline
       type &surface
       &$P_1$&$P_2$&$P_3$&$P_4$&$P_5$&$P_6$&$P_7$&$P_8$&$P_9$\\\hline\hline
        & $L_1$ & 0 & 0 & 1 & 1 & 1 & 1 & 1 & 0 & 0 \\\cline{2-11} 
      {$(2,2,2)$} 
        & $L_2$ & 1 & 1 & 0 & 0 & 1 & 1 & 0 & 1 & 0 \\\cline{2-11} 
        & $L_3$ & 1 & 1 & 1 & 1 & 0 & 0 & 0 & 0 & 1 \\\hline\hline 
        & $L_1$ & 0 & 0 & 1 & 1 & 1 & 1 & 1 & 0 & 0 \\\cline{2-11} 
        {$(2,2,4)$} 
        & $L_2$ & 1 & 1 & 0 & 0 & 1 & 1 & 0 & 1 & 0 \\\cline{2-11} 
        & $K_3$ & 1 & 1 & 1 & 1 & 0 & 1 & 1 & 1 & 2 \\\hline\hline 
    \end{tabular}
    \caption{multiplicities in the properly elliptic case}
    \label{table:multelliptic}
\end{table}

The three families of blown-up $K3$-surfaces are related via Cremona
transformations. 
The ordinary plane Cremona transformation is the rational map
defined by the linear
system of conics through three points in general position.
In suitable coordinates it can be given by the formula
$(x \cn y\cn z) \mapsto (1/x \cn 1/y\cn 1/z)$.
This formula generalises
to higher dimensions.
In particular, the space transformation,
also known as {\em reciprocal transformation},
$$
    (x\cn y\cn z \cn w)
    \mapsto \left(\frac{1}{x}\cn\frac{1}{y}\cn\frac{1}{z}
    \cn\frac{1}{w}\right)  
$$
simultaneously blows up the vertices and
blows down the faces of the coordinate tetrahedron.
The vertices are called {\em fundamental points} of
the reciprocal transformation.
Let $X\subset\PP^3$ by a surface of degree $d$
not containing any of the coordinate planes. Let $m_1$, \dots, $m_4$ be the
multiplicities of $X$ in the fundamental points.
Then the image $Y$ of $X$ is a surface of degree $3d-m_1-\cdots-m_4$.
In many cases $X$ will be singular in the fundamental points
with singularities obtained from contracting the
intersection curves of $X$ with the coordinate planes.

Specifically, a reciprocal transformation with fundamental points
$P_1$, $P_2$, $P_4$ and $P_5$ will transform a surface of type
$(4,4,4)$ into one of type $(2,4,6)$. To get from there
to a surface of type
$(2,2,8)$ we can apply a transformation with fundamental points
$P_2$, $P_5$, $P_6$ and $P_7$, where the points are as in
Table \ref{table:multiK3}.
The two other families are also related via reciprocal transformations.

\section{Families with ten triple points}
For a sextic surface $X$ with ten isolated triple points $p_g(\wt X) =0$
(\cite{EPS}, Cor.~4.6) so the ten points never lie on a quadric.
Leaving out one point the remaining nine triple points determine
a quadric $Q$.
The general element of the  pencil spanned by our  sextic and $Q^3$
is a surface with nine isolated triple points and belongs therefore
at least to one of the five families above.

\begin{lemma}\label{lem}
A sextic with ten triple points belongs to the closure of the family
of type $(2,2,2)$ or of the family of type $(4,4,4)$.
\end{lemma} 

\roep Proof.
No three triple points lie on a line (\cite{EPS}, Lemma~3.1). Two
different $(-1)$-conics meet in two triple points (\cite{EPS}, Cor.~4.8).
We study the planes containing $(-1)$-conics. 
If three planes have a line in common, there
would be $2+3\cdot3=11$ triple points; if four plane have a triple point
in common, they contain $1+6+4=11$ points, again contradicting that
the surface has ten triple points.
The number of planes is at most six. If there are exactly six,
then each triple point lies in three planes, and leaving out one of the
points
gives sextics with three planes, so of type $(2,2,2)$. 
If there are five planes, we have ten lines each containing two triple
points, so five points lie in three planes and five only in two. Leaving
out a point in only two planes gives a sextic of type $(2,2,2)$.
If there are four planes and only one point lies in three of them,
the fourth plane contains six points. So there are at least two points
in three planes each. A plane containing them both has only four
points on intersection lines so leaving out the fifth point in such a plane
gives a sextic of type $(2,2,2)$. If there are three planes, they can 
contain at most nine points, so by leaving out the point not on a plane
we keep three planes. If there are only two planes we can leave out
a point on the intersection line to get sextics without planes, so
of type $(4,4,4)$. If there is only one plane we leave out any point
in that plane.   
\qed


\boldmath
\subsection{Type $(2,2,2)$}
\unboldmath
We describe equations
for the surfaces. After a change of
coordinates we may assume that the three planes are the
sides of the coordinate tetrahedron. The remaining 
coordinate transformations are given by diagonal matrices.
We have two points on each axis in affine space and three additional
ones on the triangle at infinity. We take them to be
$P_7=(0\cn1\cn\lambda\cn0)$,
$P_8=(\mu\cn0\cn1\cn0)$ and
$P_9=(1\cn\nu\cn0\cn0)$. 
The equation has now the form
$$
    \alpha \; Q^3 + \beta \; xyztQ + \gamma \; xyz K\;,
$$
with $K$ is a four-nodal cubic passing through
$(0\cn0\cn0\cn1)$. 
With notation slightly different from \cite{EPS} we get
\begin{eqnarray*}
    Q &= & c_1c_2c_3t^2
            + t(b_1c_2c_3x+b_2c_1c_3y+b_3c_1c_2z) \\
        &&{} + c_2c_3x(\nu x - y - \mu\nu z)
          + c_1c_3y(\la y - z - \la\nu x)
          + c_1c_2z(\mu z - x - \la\mu y)\;,
\end{eqnarray*}
\begin{eqnarray*}
    K &= & t^2(\la\nu c_1x
               +\la\mu c_2y
               +\mu\nu c_3z) \\
        &&{} + t\big( \la b_1x(\nu x - y - \mu\nu z)
            +    \mu b_2y(\la y - z - \la\nu x)
            +    \nu b_3z(\mu z - x - \la\mu y)\big)\\
        &&{} +(\nu x - y - \mu\nu z)
           (\la y  - z - \la\nu x)
           (\mu z - x -\la\mu y)\;.
\end{eqnarray*}
We use the remaining freedom in coordinate transformations to
place the putative tenth triple point in $(1\cn1\cn1\cn1)$. 
We compute in the affine chart $t=1$.
The condition for a triple point is then that the function, its derivatives
and the second order derivatives vanish at $(1,1,1)$.
This gives ten equations which are linear in 
$\alpha$, $\beta $ and $\gamma$, so we may eliminate them: the
maximal minors of the coefficient matrix have to vanish.
We have
\begin{eqnarray*}
    \frac{\partial Q^3}{\partial x}           & =& 3Q^2Q_x\;,\\
    \frac{\partial^2 Q^3}{\partial x^2}        &=& 3Q^2Q_{xx}+6QQ_x^2 \;,\\
    \frac{\partial^2 Q^3}{\partial x\partial y} &=& 3Q^2Q_{xy}+6QQ_xQ_y\;.
\end{eqnarray*}
All these expressions are divisible by $Q$. 

Now we plug in $x=y=z=1$. From $Q$ we get
$$
\displaylines{
  \quad  Q(1,1,1) \hfill\cr\hfill\qquad{}=
         c_1c_2c_3  
         + c_2c_3(b_1+\nu  - 1 - \mu\nu  )
           + c_1c_3(b_2 + \la  - 1 - \la\nu  )
           + c_1c_2(b_3+ \mu  - 1 - \la\mu  ) \;,\;}
$$
an expression which we continue to denote by $Q$.
We also get  expressions for all derivatives.
Likewise we have 
\begin{eqnarray*}
    K &= & \la\nu c_1+\la\mu c_2+\mu\nu c_3 \\
        &&{}+\la b_1(\nu - 1 - \mu\nu)
     +\mu b_2(\la - 1 - \la\nu)
     +\nu b_3(\mu - 1 - \la\mu)\\
        &&{}+(\nu - 1 - \mu\nu)
      (\la - 1-  \la\nu)
      (\mu - 1 - \la\mu)\;. 
\end{eqnarray*}
Furthermore
\begin{eqnarray*}
    \left.\frac{\partial\, xyzK}{\partial x}\;\right|_{(1,1,1)}     
     &= & (yzK + xyzK_x)_{|_{(1,1,1)}}= K+K_x\;, \\
    \left.\frac{\partial^2\, xyzK}{\partial x^2}\;\right|_{(1,1,1)}
     &=& (2yzK_x + xyzK_{xx})_{|_{(1,1,1)}}= 2K_x+K_{xx} \;,\\
    \left.\frac{\partial^2 \,xyzK}{\partial x\partial y} \;\right|_{(1,1,1)} 
      & =& (zK +xzK_x+yzK_y+xyz K_{xy})_{|_{(1,1,1)}}\\
       & & \qquad\qquad {}= K +K_x+K_y+K_{xy}\;.
\end{eqnarray*}
After dividing the first row  by $Q$, which is allowed
because the tenth triple point does not lie on the quadric $Q$,
our matrix has the following form:
$$
    \pmatrix{
        Q^2 &  3QQ_x & \ldots & 3QQ_{xx}+6Q_x^2 
                                 & \ldots & 3QQ_{xy}+6Q_xQ_y & \ldots\cr
        Q   & Q +Q_x & \ldots & 2Q_x+Q_{xx}     
                                 & \ldots & Q+Q_x+Q_y+Q_{xy} & \ldots\cr
        K   & K+K_x & \ldots & 2K_x+K_{xx}     
                               & \ldots & K+K_x+K_y+K_{xy} & \ldots
    }\;.
$$
The vanishing of the maximal minors is the necessary condition
for multiplicity $3$ in the point $(1,1,1)$, but it
is not sufficient for the existence of a surface with only
isolated singularities. 
We have to cut away unwanted solutions, like $Q=Q_x=Q_y=Q_z=0$,
which makes all minors vanish, but does not give isolated triple
points.
The minors are rather formidable expressions. We first try to simplify
the matrix itself.

We start by subtracting $3Q$ times the second row from the
first row to remove all second derivatives from the first row.
After that we apply only column operations.
Some experimentation with the matrix showed that it is possible to
get two zeroes in one column. 
We observe that $Q_{xx}+2\nu Q_{xy} +\nu^2 Q_{yy}=0$. Note that 
one can write 
$Q(x,y,z,t)= \frac12 Q_{xx}x^2+Q_{xy}xy+\cdots+\frac12 Q_{tt}t^2$,
as the second derivatives are constants.
The identity $Q_{xx}+2\nu Q_{xy} +\nu^2 Q_{yy}=0$ now follows from the
fact that the point $(1\cn\nu\cn0\cn0)$ lies on the quadric.
The same point is a double point of the cubic $K$, so all
first derivatives vanish, giving by the same argument
that $(K_w)_{xx}+2\nu (K_w)_{xy} +\nu^2 (K_w)_{yy}=0$, where 
$w$ is one of $(x,y,z,t)$. Applying Euler's relation
$3K=xK_x+yK_y+zK_z+tK_t$ in the point $(1\cn 1\cn 1\cn1)$ yields by
adding that also $K_{xx}+2\nu K_{xy} +\nu^2 K_{yy}=0$, where 
now $K_{xx}$ again stands for the second derivative evaluated in
$(1,  1,1)$.
Equivalent equations
hold for the other second partials.

We get in this way three columns with two zeroes by elementary column
operations, if we  multiply one column, say the one containing
containing $Q_{xx}$, with $1+\la\mu\nu$. The vanishing of this factor
expresses  that the three points $P_7$, $P_8$ and
$P_9$ lie on a line, so we may introduce new unwanted solutions,
which we cut away later in the computation.
The result is
$$
    \pmatrix{
        -2Q^2 & -Q^2 & \ldots & Q^2-3QQ_x-3QQ_y+6Q_xQ_y & \ldots & E_\nu  
                                                          & \ldots \cr
          Q   &  Q_x & \ldots &  Q_{xy}   & \ldots & 0 & \ldots \cr
          K   &  K_x & \ldots &  K_{xy}   & \ldots & 0 & \ldots
    },
$$
where $E_\nu$ is the first  of three similar equations
$$
\begin{array}{lc}
    E_\nu\colon
       & (\nu^2+\nu+1)Q^2-3(\nu+1)Q(\nu Q_y+Q_x)+3(\nu Q_y+Q_x)^2 \;, \\
    E_\la\colon
     &(\mu^2+\mu+1)Q^2-3(\mu+1)Q(\mu Q_x+Q_z)+3(\mu Q_x+Q_z)^2 \;,\\
    E_\mu\colon
     &(\la^2+\la+1)Q^2-3(\la+1)Q(\la Q_z+Q_y)+3(\la Q_z+Q_y)^2 \;.
\end{array}
$$
These equations have to hold, for 
if $E_\nu\neq 0$, then $\alpha=0$ and the equation for the sextic is 
divisible by $xyz$. Considered as quadratic equation in $Q$ and
$\nu Q_y+Q_x$ the equation $E_\nu$ has discriminant $-3(\nu-1)^2$.
The case $\nu =1 $ is excluded:
if $\nu =1 $ then $0=Q_x+Q_y-Q=c_1c_2(c_3+b_3+\mu)$, which means that
the point $(0,0,1)$ is a triple point, which lies on the line
through the tenth point $(1,1,1)$ and $P_9=(1\cn1\cn0\cn0)$. 
Therefore no solution is defined over $\RR$. 
We have to adjoin $\sqrt{-3}$ or what amounts to the same,
the third roots of unity.


By factorising  $E_\kappa$, $\kappa=\la, \mu, \nu$, 
we get linear equations, which express
$Q_x+\nu Q_y$, $Q_y+\la Q_z$ and  $Q_z+\mu Q_x$ as  multiples of $Q$.
To express $Q_x$, $Q_y$ and $Q_z$ themselves as  multiples of $Q$
we have to multiply with the determinant  $1+\la\mu\nu$ 
of the system. By doing so to the fifth, sixth and seventh column
of our matrix we can get use the first column to get zeroes on the
first row in all other columns. This reduces our problem to
the minors of a $(2\times6)$-matrix. 

The analysis up to this point is basically contained in \cite{EPS}.
To proceed further we note that our three linear equations are in fact linear
in  $b_1c_2c_3$, $b_2c_1c_3$ and   $b_3c_1c_2$. Therefore
they can be used to eliminate the $b_i$.
For the second row this is quite easy to do:
by column operations we can remove the $b_i$ from column 5, 6 and 7
and then we take suitable linear 
combinations of columns 2, 3 and 4 with coefficients 
polynomials in $(\la,\mu,\nu)$ such that the entries on the second row
have the same coefficients at the $b_ic_ijc_k$ as our three equations. 
For the third row one has to  first multiply with a quite
complicated determinant, which leads to long expressions.
At this stage the use of the computer becomes indispensable.
The new second column turns out to be divisible by $\nu -1$,
and likewise the third by $\la-1$, the fourth by $\mu-1$.
After division the entries $(2,2)$, $(2,3)$ and $(2,4)$
are equal, which means that we again get columns with two zeroes, giving two 
equations. From the remaining $(2\times4)$-matrix we take the 6 maximal
minors. Now we have a system of 8 rather complicated
equations in 6 variables. We still have to cut away unwanted solutions,
those lying in  $Q=0$, $\la\mu\nu+1=0$, $\la=1$, $\mu=1$,
$\nu=1$ and  $c_i=0$. This can be done in {\tt Singular} as follows.
First we homogenise with an extra variable $h$. To cut away the solutions
in $Q=0$ we adjoin the inhomogeneous equation $Q-1$, where $Q$ is 
made homogeneous with $h$, and compute a standard basis. Then we homogenise
again with $Q$. By doing the same for the other  unwanted solutions
we finally obtain equations of reasonably low degree. To do the calculation
in reasonable time it is best to compute over a finite field
$\ZZ/p\ZZ$ containing the third roots of unity.
One can then try to lift the result to characteristic zero and check
whether the guessed equations really solve the system.

Let $\ep$ be a primitive third root of unity.
We first take the same root to solve the three equations $E_\kappa$:
\begin{eqnarray*}
    3(\nu Q_y+Q_x)-((1-\ep^2)\nu+(1-\ep))Q &=& 0\;, \\
    3(\mu Q_x+Q_z)-((1-\ep^2)\la+(1-\ep))Q &= &0\;,\\
    3(\la Q_z+Q_y)-((1-\ep^2)\mu+(1-\ep))Q &= &0\;.
\end{eqnarray*}
By eliminating $c_2$ and $c_3$ we end up with one equation
which is quadratic in $c_1$, so we find a three dimensional solution space.
The equations are rather involved.
%

A cyclic permutation of the variables $(x,y,z)$ in the original configuration
induces a cyclic permutation of each of the triples $(b_1,b_2,b_3)$,
$(c_1,c_2,c_3)$ and $(\la,\mu,\nu)$. A transposition of $x$ and $y$
has a more complicated effect on the coefficients.
On the points $P_7$, $P_8$, $P_9$ it acts
as  $(0\cn1\cn\lambda\cn0) \mapsto (1\cn0\cn\lambda\cn0)
=(1/\la\cn0\cn1\cn0) $,
$(\mu\cn0\cn1\cn0)\mapsto (0\cn\mu\cn1\cn0)=(0\cn1\cn1/\mu\cn0)$ and
$(1\cn\nu\cn0\cn0)\mapsto(\nu\cn1\cn0\cn0)=(1\cn1/\nu\cn0\cn0)$. 
The induced action on the coefficients is therefore
$((\la,\mu,\nu)\mapsto (1/\mu\cn1/\la\cn1/\nu)$. 
By also considering $P_1$, \dots, $P_6$ we find that 
$(c_1,c_2,c_3)\mapsto (c_2/\la\nu,c_1/\mu\nu,c_3/\la\mu)$
and $(b_1,b_2,b_3)\mapsto(b_2/\la\nu,b_1/\mu\nu,b_3/\la\mu)$.
We clear denominators in $Q(x,y,z)$ and $K(x,y,z)$.
The equation
$3(\nu Q_y+Q_x)-((1-\ep^2)\nu+(1-\ep))Q=0$ is transformed into
$3(\nu Q_y+Q_x)-((1-\ep^2)+(1-\ep)\nu)Q=0$ (where we multiplied with $\nu$
to avoid denominators). 
By taking a particular normal form of the family we found
two components, one with $\ep$ and one with $\ep^2$,
but the surfaces in those components are isomorphic.
As the permutation of $x$ and $y$ is isotopic to the identity
there is only one component (of dimension $3+15$)
in the space of all sextics.

Now we  take different roots of unity  in the  equations $E_\kappa$.
By using the permutations of $(x,y,z)$ it suffices to
consider:
\begin{eqnarray*}
    3(\nu Q_y+Q_x)-((1-\ep^2)\nu+(1-\ep))Q &=& 0\;, \\
    3(\mu Q_x+Q_z)-((1-\ep^2)\la+(1-\ep))Q &= &0\;,\\
    3(\la Q_z+Q_y)-((1-\ep)\mu+(1-\ep^2))Q &= &0\;.
\end{eqnarray*}
We start the computation  as described above. The two equations
coming from the second row of the matrix   factorise. Disregarding
a factor $(\la\mu\nu+1)$ we find the equations
$$
\displaylines{
(\la\mu c_2-c_3)\left(
\nu c_1c_2c_3-(\nu(\ep^2+\la)c_1c_3-\nu c_2c_3-\ep c_1c_2)(\mu\nu-\nu+1)
\right)\;,\cr
(\la\nu c_1-c_2)\left(
\mu c_1c_2c_3-((\ep\nu+1)c_1c_3-\ep^2\mu\nu c_2c_3-\mu c_1c_2)(\la\mu-\mu+1)
\right)\;.
}
$$
Applying a suitable transposition of the coordinates induces a
transformation which sends the first equation to the second one
with $\ep$ replaced by $\ep^2$.
 
We find one three dimensional solution  by taking both long factors.
Then we find  $\mu c_2=(\la\mu-\mu+1)(\mu\nu-\nu+1)$
and  a quadratic equation in $c_1$, which I do not describe here.

Another three dimensional solution is found by taking the equations
$\la\mu c_2-c_3$ and 
$\mu c_1c_2c_3-((\ep\nu+1)c_1c_3-\ep^2\mu\nu c_2c_3-\mu c_1c_2)(\la\mu-\mu+1)$.
The other possible choice gives a solution, isomorphic to the
complex conjugate of this one. It might seem that we get 
two different solutions, but as we shall show, the surfaces in 
question can also be written in a different way as
a degeneration of a sextic of type $(2,2,2)$. 
By computing for a specific example we find that both solutions
are slices of the same component in the space of all sextics.
This time we find a linear equation for $c_1$:
$$
c_1 +\ep^2(\la\mu-\mu+1)(\la\mu\nu+\ep\mu\nu-\ep\nu-\ep^2)=0\;.
$$
We have already $c_3= \la\mu c_2$ and we find
$$
(\la\mu\nu+\ep\mu\nu+\ep^2\nu-\ep^2)c_3+\ep(\ep \la\nu+\la-1)c_1=0\;.
$$

Finally, taking $c_3=\la\mu c_2$ and $c_2=\la\nu c_1$ gives a 
two dimensional solution consisting of two components, one of which lies
inside the last component just found, and the other in the one obtained 
by interchanging the equations.

\begin{proposition}
The  family of sextic surfaces  of type
$(2,2,2)$ with nine triple points contains in its closure 
three different families
of sextics with ten triple points, which contain three, four or five
$(-1)$-conics.
\end{proposition}

\roep Proof.
We have already seen that there are at most three different families.
We distinguish between them with the number of $(-1)$-conics.
A $(-1)$-conic determines a plane, whose intersection with one of the three
coordinate planes can have at most two triple points.
It has to contain at least two  of the points $P_1$, \dots, $P_6$ 
on the coordinate axes,  because
$P_7$, $P_8$, $P_9$ and $P_{10}$ are not coplanar. But if the plane contains
a point on a coordinate axis, it contains only two other triple points on the 
coordinate planes through the point and therefore it contains the two
points not in these planes.  If there are three points of the points
$P_1$, \dots, $P_6$ in the plane, it therefore contains again  $P_7$,
$P_8$, $P_9$ and $P_{10}$. 
Therefore there are only three possible planes
which can contain a $(-1)$-conic, namely the planes through $P_{10}$
and two of $P_7$, $P_8$ and $P_9$.
The equation for the plane through $P_7$, $P_8$  and $P_{10}$ is
$\mu z-x-\la \mu y+\la\mu-\mu+1$.

To determine the number of $(-1)$-conics in each family it suffices to do it
for a specific example. We obtain three conditions by requiring
that the points $(1\cn0\cn0\cn1)$,  $(0\cn1\cn0\cn1)$ and $(0\cn0\cn1\cn1)$
are triple points. This gives the equations
$c_1+b_1+\nu=0$, $c_2+b_2+\la=0$ and $c_3+b_3+\mu=0$.
In the first family we find $\la=\mu=\nu$, $c_1=c_2=c_3$,
$\nu^4-3\ep\nu^2+\ep^2=(\nu^2+\ep^2\nu-\ep)(\nu^2-\ep^2\nu-\ep)$, 
$c_3^2+(1-\ep^2)\nu^2-\ep+1$.
In the last family found above we get
$\la=\nu$, $\mu-2\ep^2\nu+3$, $c_1=c_3$, $c_3+(\ep^2-1)\nu+\ep-1$,
$\nu^2-\ep\nu-\ep^2$ and  $c_2+(\ep-\ep^2)\nu$.
For the third family this specialisation does not work, so we have
to take a different one.
By checking in finite characteristic we make sure that there really 
exist a sextic with ten isolated triple points for these parameter values.

We then determine if one of the three planes contains more than three
triple points. The result is that 
the first family does not contain extra $(-1)$-conics.
The second family contains one extra $(-1)$-conic, the plane through
$P_7$, $P_8$  and $P_{10}$, which also contains a point on the $x$-axis
and on the $y$-axis, with coordinates $(c_1\cn 0\cn 0\cn \nu)$
resp.~$(0\cn c_2 \cn0\cn l)$.

We specialise the third family by  taking suitable values for $\la$,
$\mu$ and $\nu$. A good choice is $\la=\mu=-1$, $\nu=\ep$.
We can then compute the intersection points of the three planes with the
coordinate axes and check whether they are triple points.
We find the equations $3b_1+c_1+9\ep$, $b_2+2\ep-4$, $c_2-6\ep+3$,
$(\ep+3)b_3+c_3-5\ep-8$ and a
quadratic equation for $c_1$, which does not factor in an easy way.
For both values of $c_1$ 
the two points on the $y$-axis are given by
$(y-3)(y+2\ep-1)$, on the $x$-axis lies $(3,0,0)$ and on the
$z$-axis $(0,0,\ep+3)$. The result is that there are two extra $(-1)$-conics,
the one through $P_7$, $P_8$  and $P_{10}$
and the one through $P_8$, $P_9$  and $P_{10}$.
\qed

\roep Remark.
The computation shows that there are no sextics with ten isolated
triple points and six $(-1)$-conics. The arguments proving  Lemma \ref{lem}
do not exclude such a configuration. In fact we can take the three
planes in the proof above and take as the points on the coordinate
axes the intersection points with these planes. But a sextic with
these isolated triple points occurs in a pencil, containing also the
product of the six planes. The matrix above should then have rank one.
The first $2\times2$ minor gives the equation $Q^2(Q-2Q_x)=0$,
so together with the equations $E_\kappa$ we find $Q=0$, contradicting
the fact that the ten points do not lie on a quadric.

\boldmath
\subsection{Type $(4,4,4)$}
\unboldmath

To complete the classification of sextics with ten triple points
we look for a tenth triple point in the family of type $(4,4,4)$.
Equations for the family are given in \cite{EPS}, which depend on seven
parameters. It is convenient to work with more parameters,
which then allows to take the tenth point in fixed position.

We take three quadratic cones $K_i$
with vertices $P_{i+1}$ at infinity such that
$K_i$ passes through $P_{i-1}$ but not through
$P_{i}$, where we compute the indices modulo 3. 
In general the quadrics intersect in eight distinct points.
We require that two of them are the points $(0,0,0)$ 
and $(u, v,w)$. The six remaining points will be 
the triple points of the sextic. We find
\begin{eqnarray*}
    K_1 &=& wx^2+auz+bwx-(a+b+u)xz\;,\\
    K_2 &=& uy^2+cvx+duy-(c+d+v)xy\;,\\
    K_3 &=& vz^2+ewy+fvz-(e+f+w)yz\;.
\end{eqnarray*}
To compute $Q$, the quadric through $P_1$, \dots, $P_9$, 
but not through $(0,0,0)$ and $(u, v, w)$,
we note that
the $K_i$ lie in the ideal $(u-x,v-y,w-z)$. We can write
$$
    \pmatrix{
        K_3 \cr
        K_1  \cr
        K_2  
    } =
    \pmatrix{
        0        & (f+z)z  &  (e-z)y \cr
        (a-x)z   & 0       &  (b+x)x \cr
        (d+y)y   & (c-y)x  &  0 
    }
    \pmatrix{
        u-x \cr
        v-y  \cr
        w-z  
    }
    \;.
$$
Dividing the determinant of the matrix by $xyz$ gives the
inhomogeneous equation 
$$
    Q=(a-x)(c-y)(e-z)+(b+x)(d+y)(f+z)
$$
which is indeed the sought quadric.
%
Note that our equations are  homogeneous in the coefficients
$a$, \dots, $w$ and the affine coordinates $x$, $y$, $z$ together.

The obvious thing to do now is to determine the conditions under which
a surface $\la K_1K_2K_3 + \mu Q^3$ has a triple point in
$(x,y,z)=(1,1,1)$. Despite great efforts I did not succeed in finding
a single example. Finally I decided to compute the transformations
which bring the known example from \cite{EPS} (which is the same as
the specific example in the first family above) into this family.
The result was that the tenth point lies in the plane at infinity.
In fact, a  long, but doable computation with
{\tt Singular} shows that only solutions of the equations 
occur when $(1,1,1)$ lies on the quadric $Q$ or one of the cones $K_i$.

We therefore now search  under the 

\roep Assumption.
The point $(1\cn1\cn1\cn0)$ is a triple point.
\endroep

For the pencil $\alpha\, K_1K_2K_3 + \beta\, Q^3$ we compute all ten
second partial derivatives and evaluate them in $(1\cn1\cn1\cn0)$. 
The resulting equations are linear in $\alpha$ and $\beta$,
so we  eliminate these variables and end up with a $2\times10$ matrix.

The vanishing of the minors of the matrix is again a necessary
condition, for the existence of a sextic with ten triple points,
but it is not sufficient for 
isolated triple points.
Indeed, there are some easy to see `false' solutions:
if $K_1=K_2=K_3=0$, then the whole first row vanishes (we take
at most second derivatives of the product $K_1K_2K_3$) and we get 
$\beta=0$. 
Also, if $a+b=c+d=e+f=0$, the second row vanishes. We know that 
the ten points cannot lie on the quadric $Q$.
We only want solutions with $Q\neq0$, 
$K_1\neq0$, $K_2\neq0$ and $K_3\neq0$.

Our equations are homogeneous in $a$, \dots, $z$. Moreover, the 
derivatives not involving $t$ depend only on the sums $a+b$, $c+d$, $e+f$
and the $u$, $v$, $w$: note that $Q|_{t=0}=(a+b)yz+(c+d)xz+(e+f)xy$ and
$K_1|_{t=0}=wx^2-(a+b+u)xz$. 
This means that we can start by analysing the six first columns. 
We cut away one after another the solutions lying in  $Q=a+b+c+d+e+f=0$,
$K_1=w-(a+b+u)=0$, $K_2=0$ and $K_3=0$.
To dispose of the solutions in a hyperplane $L=0$ we add the inhomogeneous 
equation $L=1$ and compute a standard basis. Afterwards we make the
equations homogeneous again.

The computation with {\tt Singular}  gives twelve equations. 
They define two complex conjugate 
components. Eliminating $a$, $c$ and $e$ gives two equations
$$
\displaylines{
(b+d+f)^2+(b+d+f)(u+v+w)+(u+v+w)^2\;,\cr
uv+uw+vw\;.}
$$
Again we have to adjoin the third roots of unity. 
With $\ep$ a primitive third
root of unity  we find two components,
one of them given by
$$
\displaylines{
e+f+\ep^2v-\ep w\cr
c+d+\ep^2u-\ep v\cr
a+b+\ep^2w-\ep u\cr
b+d+f-\ep(u+v+w)\cr
uv+uw+vw\;.}
$$

We give an explicit example: $v=w=2$, $u=d=-1$, $b=0$, so
\begin{eqnarray*}
    K_1 &=& 2x^2-(\ep+2)z+\ep^2xz\;,\\
    K_2 &=& -y^2+2\ep x+y+\ep^2 xy\;,\\
    K_3 &=& 2z^2-2\ep^2 y+(6\ep+2)z+4\ep^2 yz\;,\\
    Q   &=& -(\ep+2-x)(\ep-y)(\ep^2+z)+x(y-1)(z+3\ep+1)\;.
\end{eqnarray*}
Then  $27\, K_1K_2K_3 + 2\, Q^3$ has ten ordinary triple points.
To find them it is convenient to compute in finite characteristic $p$.
After some experimentation I found that for $p=67$ with $\ep=-30$ 
all points are defined over the base field.

\begin{proposition}
The  family of sextic surfaces  of type $(4,4,4)$  contains in its closure 
one family
of sextics with ten triple points, which each contain two
$(-1)$-conics.
\end{proposition}

\roep Proof.
One of the intersection points of the quadric cones $K_i$ lies
in the plane $t=0$. To see this
we observe that $K_3|_{t=0}= z(vz  + \ep^2(v+w)y)$. By cyclic
permutation we get three lines $vz  + \ep^2(v+w)y$,
$wx  + \ep^2(w+u)z$ and $uy  + \ep^2(u+v)z$. The condition that they
pass though one point is 
$$ 
(u+v)(v+w)(w+u)+ uvw = (u+v+w)(uv+uw+vw)=0\;,
$$
which is satisfied on our component. 

The intersection point
is $(\ep^2 uw \cn \ep uv \cn vw \cn 0)$. 
Together with the tenth point $(1\cn1\cn1\cn0)$ it lies on the 
line $t=x+\ep y+\ep^2 z=0$. One of the planes in the pencil of planes through
this line contains three more triple points. It can be found
by transforming the coordinates $(x,y,z)$  into the eigenfunctions
of cyclic permutation, making  $x+\ep y+\ep^2 z$ into a coordinate
and eliminating the others. The computation is best done
in finite characteristic.
Once the result is known one can find a derivation.
We observe the following factorisation
modulo the ideal defining the component
$$
uK_3+\ep^2 vK_1+\ep wK_2\equiv 
       ( \ep^2 vwx+wuy+\ep uvz)(x+\ep y+\ep^2z+e+\ep d)\;.
$$
In the affine chart $t\neq 0$ the six common points of the
quadric cones lie therefore on two planes. The first factor
contains the point $(u, v , w)$, while the second factor
is the plane of the pencil we are after. Note also that 
$x+\ep y+\ep^2z+e+\ep d$ and $y+\ep z+\ep^2x+a+\ep f$
give the same plane.

By leaving out the point $P_9= (\ep^2 uw \cn \ep uv \cn vw \cn 0)$
we realise our surface in a different way as special element in a pencil
of type $(4,4,4)$. A coordinate transformation
brings it in our standard form. To determine it we have to know the position 
of the three vertices, so we only compute in our specific example.
We obtain values for the parameters $(a,\dots,f,u,v,w)$ and 
compute that they satisfy the equations for the complex conjugate
component. This shows that there is only one family. 
\qed

\subsection{Cremona transformations}
To compute  the effect of a Cremona transformation it is useful
to know about other $(-1)$-curves on our surfaces.
Each family lies also in the closure of  other families of sextics
with nine triple points. For explicit computations we
need to know the coordinates of the ten triple points. We use the
specific examples in finite characteristic. 

We start with the surface with two $(-1)$-conics.
If we leave out $P_1$, then the surfaces
has one $(-1)$-conic, so is of
of type $(2,4,6)$ with the $(-1)$-conic the one determined above.
The pencil has to contain the reducible surface $L_2K_1C_1$ with
$C_1$ a cubic surface. In the example one finds an explicit equation
for $C_1$. 
Leaving out $P_7$ or $P_8$ gives a surface with two $(-1)$-conics,
which a priori can be of type $(2,2,8)$ or $(2,2,4)$. The explicit
example shows that the first case occurs.
Table \ref{table:multten} contains all the surfaces found in this way, with
$L_i$ planes, $K_i$ quadric cones, $C_i$ four-nodal cubics and $Q_i$
quartics with one triple point and six nodes.
Through each point pass 13 of the 16 surfaces and the reducible surface in 
the pencil obtained by leaving out this point is the union of the
other three surfaces.

  
\begin{table}[h]
    \centering
    \begin{tabular}{|c||c|c|c|c|c|c|c|c|c|c|}\hline
        surface
       &$P_1$&$P_2$&$P_3$&$P_4$&$P_5$&$P_6$&$P_7$&$P_8$&$P_9$&$P_{10}$
      \\\hline\hline
        $L_1$ & 1 & 1 & 1 & 0 & 0 & 0 & 0 & 0 & 1 & 1  \\\hline
        $L_2$ & 0 & 0 & 0 & 1 & 1 & 1 & 0 & 0 & 1 & 1  \\\hline\hline
        $K_1$ & 0 & 2 & 1 & 1 & 1 & 1 & 1 & 1 & 1 & 0  \\\hline
        $K_2$ & 1 & 0 & 2 & 1 & 1 & 1 & 1 & 1 & 1 & 0  \\\hline
        $K_3$ & 2 & 1 & 0 & 1 & 1 & 1 & 1 & 1 & 1 & 0  \\\hline
        $K_4$ & 1 & 1 & 1 & 0 & 2 & 1 & 1 & 1 & 0 & 1  \\\hline
        $K_5$ & 1 & 1 & 1 & 1 & 0 & 2 & 1 & 1 & 0 & 1  \\\hline
        $K_6$ & 1 & 1 & 1 & 2 & 1 & 0 & 1 & 1 & 0 & 1  \\\hline\hline
        $C_1$ & 0 & 1 & 2 & 1 & 1 & 1 & 2 & 2 & 1 & 2  \\\hline
        $C_2$ & 2 & 0 & 1 & 1 & 1 & 1 & 2 & 2 & 1 & 2  \\\hline
        $C_3$ & 1 & 2 & 0 & 1 & 1 & 1 & 2 & 2 & 1 & 2  \\\hline
        $C_4$ & 1 & 1 & 1 & 0 & 1 & 2 & 2 & 2 & 2 & 1  \\\hline
        $C_5$ & 1 & 1 & 1 & 2 & 0 & 1 & 2 & 2 & 2 & 1  \\\hline
        $C_6$ & 1 & 1 & 1 & 1 & 2 & 0 & 2 & 2 & 2 & 1  \\\hline\hline
        $Q_1$ & 2 & 2 & 2 & 2 & 2 & 2 & 0 & 3 & 1 & 1  \\\hline
        $Q_2$ & 2 & 2 & 2 & 2 & 2 & 2 & 3 & 0 & 1 & 1   \\\hline\hline
    \end{tabular}
    \caption{multiplicities of the $(-1)$-curves at the singular points}
    \label{table:multten}
\end{table}

To get  with a Cremona transformation
again a surface with ten isolated triple points
we have to take the four fundamental points such that no three lie in
a plane. 
For the surfaces of type $(4,4,4)$ there are only a few possibilities,
due to the symmetry in the configuration. 
We can compute the strict transform of
each of the surfaces in Table \ref{table:multten} using the degree
formula $3d-m_1-\cdots-m_4$. The multiplicity of the transformed surface
in one of the four image points
is the degree of the exceptional curve, which is itself the image
under a standard plane Cremona transformation of the intersection
curve of the surface with the plane through the three opposite
fundamental points: the new multiplicity $m_1$ is $2d-m_2-m_3-m_4$.

If we take $P_1$, $P_7$, $P_8$ and $P_9$ as fundamental points
the plane $L_1$ is transformed in a plane, as is the quadric $K_3$.
We get again a sextic with two $(-1)$-conics. The transform of each
of the cubics $C_4$, $C_5$, $C_6$ is a quadric cone not passing through
the new $P_1$ and simply through $P_7$, $P_8$ and $P_9$. So leaving
out the new $P_1$ gives a surface of type $(4,4,4)$ again.

We get three $(-1)$-conics if we take $P_1$, $P_2$, $P_4$ and $P_7$
as fundamental points. 
For four $(-1)$-conics  can we take $P_1$, $P_2$, $P_4$ and $P_5$
as fundamental points. 

\begin{table}[h]
    \centering
    \begin{tabular}{|c||c|c|c|c|c|c|c|c|c|c|}\hline
        surface
       &$P_1$&$P_2$&$P_3$&$P_4$&$P_5$&$P_6$&$P_7$&$P_8$&$P_9$&$P_{10}$
      \\\hline\hline
$L_1$ & 0 & 0 & 1 & 1 & 1 & 1 & 1 & 0 & 0 & 0  \\\hline
$L_2$ & 1 & 1 & 0 & 0 & 1 & 1 & 0 & 1 & 0 & 0  \\\hline
$L_3$ & 1 & 1 & 1 & 1 & 0 & 0 & 0 & 0 & 1 & 0  \\\hline
$L_4$ & 0 & 0 & 0 & 1 & 0 & 1 & 0 & 1 & 1 & 1  \\\hline
$L_5$ & 0 & 1 & 1 & 0 & 0 & 0 & 1 & 1 & 0 & 1  \\\hline\hline
$Q_1$ & 2 & 2 & 1 & 1 & 2 & 2 & 2 & 0 & 2 & 3  \\\hline
$Q_2$ & 2 & 1 & 1 & 2 & 3 & 0 & 2 & 2 & 2 & 2  \\\hline
$Q_3$ & 2 & 1 & 2 & 0 & 2 & 2 & 2 & 1 & 3 & 2  \\\hline
$Q_4$ & 2 & 2 & 0 & 2 & 2 & 1 & 3 & 1 & 2 & 2  \\\hline
$Q_5$ & 3 & 0 & 2 & 1 & 2 & 1 & 2 & 2 & 2 & 2  \\\hline\hline
    \end{tabular}
    \caption{multiplicities of the $(-1)$-curves in case of 5 planes}
    \label{table:tenfive}
\end{table}
A surface with five $(-1)$-conics cannot be obtained directly with a
reciprocal transformation. Instead we first study the configuration
in more detail. Leaving out a point on two planes gives again surfaces
of type $(2,2,2)$, whereas leaving out one of the five points 
on three planes leads to sextics of type $(2,2,8)$. 
There are five quartic surfaces $Q_i$ with a triple point.
Table \ref{table:tenfive} gives the
multiplicities of the surfaces involved at the singular points.

A Cremona transformation with fundamental points $P_1$, $P_5$, $P_7$
and $P_9$ brings us to the family with three $(-1)$-conics.
This shows that all four families are related by 
Cremona transformations (obtained by composition of reciprocal
transformations).

\vfill
\noindent Matematik\\
  Chalmers tekniska h\"ogskola och G\"oteborgs universitet, \\
  SE 412 96 G\"oteborg, Sweden\\
  e-mail: \texttt{stevens@math.chalmers.se}


\begin{thebibliography}{EPS}


\bibitem[En]{surf}
Stefan Endra\ss\ et al.,
      \emph{surf} 1.0.3 - visualizing algebraic curves and algebraic surfaces
       (2001). {\tt http://surf.sourceforge.net/}.

\bibitem[EPS]{EPS}
Stefan Endra\ss,  Ulf Persson and Jan Stevens.
       {\sl  Surfaces with Triple Points.\/} 
       J. Algebraic Geom. {\bf12}
       (2003), 307--320. 

\bibitem[GPS]{Singular}
G.-M. Greuel, G. Pfister  H. Sch\"onemann.
      \emph{Singular} 2.0. A Computer Algebra System for Polynomial
       Computations. Centre for Computer Algebra, University of
       Kaiserslautern (2001).\\ {\tt http://www.singular.uni-kl.de}.

\end{thebibliography}
\end{document}